\documentclass[10pt,fleqn,leqno]{article}
\usepackage{amsmath}
\usepackage{amssymb}
\newtheorem{theorem}                 {Theorem}      [section]
\newtheorem{proposition}  [theorem]  {Proposition}

\newtheorem{lemma}        [theorem]  {Lemma}

\title{Other Classes of Tangent Bundles with General Natural Almost
Anti-Hermitian Structures}
\author{S.~L.~Dru\c t\u a}

\date{}

\begin{document}
\maketitle

\begin{abstract} \normalsize\bf Abstract.
\rm We continue the study initiated by Oproiu in
\cite{OproiuAntiHerm}, concerning the anti-Hermitian structures of
general natural lift type on the tangent bundles. We get the
conditions under which these structures are in the eight classes
obtained by Ganchev and Borisov in \cite{Ganchev}. We complete the
characterization of the general natural anti-K\"ahlerian
structures on the tangent bundles with necessary and sufficient
conditions, then we present some results concerning the remaining
classes.
\end{abstract}

\vskip3mm

{\bf M. S. C. 2000:} primary 53C05, 53C15, 53C55.


\numberwithin{equation}{section}

\section{Introduction}

In the last fifty years, a lot of papers were dedicated to the
geometric structures obtained by lifting the metric from the base
manifold to the tangent bundle. The first Riemannian metric on the
tangent bundle was constructed by Sasaki in \cite{Sasaki}, but in
the most cases, the study of some geometric properties of the
tangent bundle endowed with this metric led to the flatness of the
base manifold. In the next years, the authors were interested in
finding of other lifted structures on the tangent bundles, with
quite interesting properties (see \cite{Anastasiei},
\cite{Cruceanu}, \cite{IanusUdriste}, \cite{Janyska},
\cite{KowalskiSek}, \cite{Oproiu10}, \cite{Udriste} -
\cite{YanoIsh}).

The results concerning the natural lifts (see \cite{Janyska} and
\cite{KowalskiSek}), allowed Oproiu to introduce on the tangent
bundle $TM$, a natural almost complex structure $J$ and a natural
metric $G$, both of them being obtained as diagonal lifts of the
Riemannian metric $g$ from the base manifold (see \cite{Oproiu3}).
The same author generalized these lifts in  \cite{Oproiu},
introducing the notion of general natural lift on the tangent
bundle.

In several recent works, like \cite{BejanOpr} - \cite{Eni},
\cite{Gartu}, \cite{Labbi}, \cite{MunteanuM} - \cite{Munteanu1},
there are studied some new geometric structures on the tangent
bundle $TM$ of a Riemannian manifold $(M,g)$, obtained by
considering the natural lifts of $g$ to $TM$.

In \cite{OproiuAntiHerm}, Oproiu gave the characterization for the
anti-Hermitian structures of ge-\ neral natural lift type on the
tangent bundle and obtained some necessary conditions under which
these structures are anti-K\"ahlerian.

In the present paper we give the complete characterization of the
anti-K\"ahlerian structures of general natural lift type on the
tangent bundle and we obtain the conditions under which the
general natural almost anti-Hermitian tangent bundles from
\cite{OproiuAntiHerm} are in the other classes of anti-Hermitian
manifolds (almost complex manifolds with Norden metric),
determined  by Ganchev and Borisov in \cite{Ganchev}.
Particularizing the obtained results to the case of the natural
diagonal or anti-diagonal lifted structures, we get the examples
constructed by Oproiu and Papaghiuc in \cite{OproiuPapII} and
\cite{OproiuPapI}.

The same authors studied in \cite{OproiuPap} the Einstein property
for the quasi-anti-K\"ahler manifolds from \cite{OproiuPapII}, and
in \cite{Papag} Papaghiuc treated a particular case of
\cite{OproiuPapI}.

The manifolds, tensor fields and other geometric objects
considered in this paper are assumed to be differentiable of class
$C^\infty$ (i.e. smooth). The Einstein summation convention is
used throughout this paper, the range of the indices
$h,i,j,k,l,m,r, $ being always $\{1,\dots ,n\}$.

\section{Preliminary results}
Let $(M,g)$ be a smooth $n$-dimensional Riemannian manifold and
denote its tangent bundle by $\tau :TM\rightarrow M$. The total
space $TM$ has a structure of a $2n$-dimensional smooth manifold,
induced from the smooth manifold structure of $M$. This structure
is obtained by using local charts on $TM$ induced  from the usual
local charts on $M$. If $(U,\varphi )= (U,x^1,\dots ,x^n)$ is a
local chart on $M$, then the corresponding induced local chart on
$TM$ is $(\tau ^{-1}(U),\Phi )=(\tau ^{-1}(U),x^1,\dots , x^n,$
$y^1,\dots ,y^n)$, where the local coordinates $x^i,y^j,\
i,j=1,\dots ,n$, are defined as follows. The first $n$ local
coordinates of a tangent vector $y\in \tau ^{-1}(U)$ are the local
coordinates in the local chart $(U,\varphi)$ of its base point,
i.e. $x^i=x^i\circ \tau$, by an abuse of notation. The last $n$
local coordinates $y^j,\ j=1,\dots ,n$, of $y\in \tau ^{-1}(U)$
are the vector space coordinates of $y$ with respect to the
natural basis in $T_{\tau(y)}M$ defined by the local chart
$(U,\varphi )$. Due to this special structure of differentiable
manifold for $TM$, it is possible to introduce the concept of
$M$-\emph{tensor field} on it (see \cite{Mok}), called by R. Miron
and his collaborators \emph{distinguished tensor field} or
$d$-\emph{tensor field}. The algebra of $d$-tensor fields on the
tangent bundle of a manifold is studied in \cite{MironAnastas}
(see also \cite{Bucataru}, \cite{MAB}).

Denote by $\dot \nabla$ the Levi Civita connection of the
Riemannian metric $g$ on $M$. Then we have the direct sum
decomposition
\begin{equation}\label{splitting}
~~~~~~~~~~~~~~~~~~~~~~~~~~~~~~~~~TTM=VTM\oplus HTM
\end{equation}
of the tangent bundle to $TM$ into the vertical distribution
$VTM={\rm Ker}\ \tau_*$ and the horizontal distribution $HTM$
defined by $\dot \nabla $ (see \cite{YanoIsh}). The set of vector
fields $\{\frac{\partial}{\partial y^1}, \dots ,
\frac{\partial}{\partial y^n}\}$ on $\tau ^{-1}(U)$ defines a
local frame field for $VTM$ and for $HTM$ we have the local frame
field $\{\frac{\delta}{\delta x^1},\dots ,\frac{\delta}{\delta
x^n}\}$, where $ \frac{\delta}{\delta
x^i}=\frac{\partial}{\partial x^i}-\Gamma^h_{0i}
\frac{\partial}{\partial y^h},\ \ \ \Gamma ^h_{0i}=y^k\Gamma
^h_{ki},$ and $\Gamma ^h_{ki}(x)$ are the Christoffel symbols of
$g$.

The set $\{\frac{\partial}{\partial y^i},\frac{\delta}{\delta
x^j}\}_{i,j=\overline{1,n}}$ defines a local frame on $TM$,
adapted to the direct sum decomposition (\ref{splitting}).

Consider the energy density of the tangent vector $y$ with respect
to the Riemannian metric $g$
\begin{equation}
~~~~~~~~~~~~~~~t=\frac{1}{2}\|y\|^2=\frac{1}{2}g_{\tau(y)}(y,y)=\frac{1}{2}g_{ik}(x)y^iy^k,
\ \ \ y\in \tau^{-1}(U).
\end{equation}
Obviously, we have $t\in [0,\infty)$ for all $y\in TM$.

We shall use the following lemma, which may be proved easily.

\begin{lemma}\label{lema1}
If $n>1$ and $u,v$ are smooth functions on $TM$ such that
$$
u g_{ij}+v g_{0i}g_{0j}=0,
$$
on the domain of any induced local chart on $TM$, then $u=0,\
v=0$. We used the notation $g_{0i}=y^hg_{hi}$.
\end{lemma}

Oproiu introduced in \cite{Oproiu} a natural $1$-st order almost
complex structure on $TM$, just like the natural $1$-st order
lifts of $g$ to $TM$ are obtained in \cite{KowalskiSek}:
\begin{eqnarray}\label{sist4}
~~~\begin{array}{c} JX^H_y=a_1(t)X^V_y+b_1(t)g_{\tau(y)}(X,y)
y^V_y+a_4(t)X^H_y+b_4(t)g_{\tau(y)}(X,y)y^H_y,
\\\newline
JX^V_y=a_3(t)X^V_y+b_3(t)g_{\tau(y)}(X,y)
y^V_y-a_2(t)X^H_y-b_2(t)g_{\tau(y)}(X,y)y^H_y,
\end{array}
\end{eqnarray}
$\forall X\in \mathcal{T}^1_0(TM),\ \forall y\in TM,$ $a_\alpha,\
b_\alpha,\ \alpha=\overline{1,4}$ being smooth functions of the
energy density, t.

In the mentioned paper, Oproiu proved the following results.
\begin{theorem}\label{th4}
{The natural tensor field $J$ of type $(1,1)$ on $TM$, given by
\rm(\ref{sist4}) \it defines an almost complex structure on $TM$,
if and only if $a_4=-a_3,\ b_4=-b_3,$ and the coefficients $a_1,\
a_2,\ a_3,\ b_1,\ b_2$ and $b_3$ are related by
\begin{equation}\label{inlocuire}
~~~~~~~~~~a_1a_2=1+a_3^2\ ,\ \ \
(a_1+2tb_1)(a_2+2tb_2)=1+(a_3+2tb_3)^2.
\end{equation}}
\end{theorem}

\begin{theorem}\label{th3}
Let $(M,g)$ be an $n(>2)$-dimensional connected Riemannian
manifold. The almost complex structure $J$ defined by
{\rm(\ref{sist4})} on $TM$ is integrable if and only if $(M,g)$
has constant sectional curvature $c$, and the  coefficients $b_1,\
b_2,\ b_3$ are given by:
\begin{equation}\label{integrab}
~~~~~\begin{array}{c} b_1=\frac{2 c^2 t a_2^2+2 c t a_1
a_2^\prime+a_1 a_1^\prime -c+3 c a_3^2}{a_1-2 t a_1^\prime-2 c t
a_2-4 c t^2 a_2^\prime},\ b_2=\frac{2 t a_3^{\prime 2}- 2 t
a_1^\prime a_2^\prime+c a_2^2+2 c t a_2 a_2^\prime+a_1
a_2^\prime}{a_1-2 t a_1^\prime-2 c t a_2-4
c t^2 a_2^\prime},\\
b_3 =\frac{a_1 a_ 3^\prime+ 2 c a_2 a_3+4 c t a_2^\prime a_3-
 2 c t a_ 2 a_3^\prime}{a_1-2 t a_1^\prime-2 c t a_2-4 c t^2
 a_2^\prime}.
\end{array}
\end{equation}
\end{theorem}

Let $G$ be the $1$-st order natural semi-Riemannian metric $G$, of
signature $(n,n)$ on $TM$, considered by Oproiu in
\cite{OproiuAntiHerm}:
\begin{eqnarray}\label{defG}
~~~~~~~~\begin{array}{c} G(X^H_y,Y^H_y)=c_1(t)g_{\tau(y)}(X,Y)+
d_1(t)g_{\tau(y)}(X,y)g_{\tau(y)}(Y,y),
\\
G(X^V_y,Y^V_y)=
c_2(t)g_{\tau(y)}(X,Y)+d_2(t)g_{\tau(y)}(X,y)g_{\tau(y)}(Y,y),
\\
G(X^V_y,Y^H_y)=c_3(t)g_{\tau(y)}(X,Y)+d_3(t)g_{\tau(y)}(X,y)g_{\tau(y)}(Y,y).
\end{array}
\end{eqnarray}
$\forall X,Y\in \mathcal{T}^1_0(TM),\ \forall y\in TM,$ where
$c_\alpha,\ d_\alpha,\ \alpha=\overline{1,3}$ are six smooth
functions of the energy density on $TM$.

The conditions for $G$ to be nondegenerate are assured if
$$
c_1c_2-c_3^2\neq0, \ (c_1+2td_1)(c_2+2td_2)-(c_3+2td_3)^2\neq 0.
$$

The almost anti-Hermitian or Norden metric with respect to the
almost complex structure $J$, was defined in \cite{Ganchev} and
\cite{OproiuPapI} as the semi-Riemannian metric $G$ satisfying
\begin{equation}\label{defantiherm}
~~~~~~~~~~~~~~~~~~~G(JX,JY)=-G(X,Y),\ \forall X,Y\in
\mathcal{T}^1_0(TM).
\end{equation}

The study of the conditions under which the metric $G$ is almost
anti-Hermitian with respect to the almost complex structure $J$,
led in \cite{OproiuAntiHerm} to the result:

\begin{theorem}\label{antiherm}
The general natural semi-Riemannian metric $G$ is almost
anti-Hermi-\ tian with respect to the general natural almost
complex structure $J$ on the tangent bundle $TM$, if and only if
the coefficients of $J$ and $G$ satisfy
\begin{eqnarray}\label{antiherm1}
~~~~~~~~~~~~~~~~~~~~~~~~~~~~~~~~~a_2c_1+a_1c_2=2 a_3c_3,
\end{eqnarray}
\begin{eqnarray}\label{antiherm2}
~~~~~~(a_2+2tb_2)(c_1+2td_1)+(a_1+2tb_1)(c_2+2td_2)=2(a_3+2tb_3)(c_3+2td_3).
\end{eqnarray}
\end{theorem}

In the same paper, Oproiu initiated the study of the anti-K\"ahler
structures of general natural lift type on $TM$, obtaining some
necessary conditions under which the almost anti-Hermitian
structures characterized by theorem \ref{antiherm} are
anti-K\"ahler.

Since in the case of the almost anti-Hermitian structures of
natural diagonal lift type on $TM$ the coefficients $a_3,\ b_3,\
c_3,\ d_3$ vanish, the relations (\ref{antiherm1}) and
(\ref{antiherm2}) become simpler, allowing the complete
characterization of the natural diagonal anti-K\"ahler structures
on $TM$ (see \cite{OproiuAntiHerm}), and the characterizations of
all the other classes of almost anti-Hermitian structures of
natural diagonal lift type on $TM$ (see \cite{OproiuPapII}).

In the following sections, we shall get some necessary and
sufficient conditions under which the general natural almost
anti-Hermitian structures on $TM$, obtained in
\cite{OproiuAntiHerm}, are in the eight classes determined by
Ganchev and Borisov in \cite{Ganchev}.

\section{Anti-K\"ahler structures of general natural lift type on $TM$}

The first class of almost anti-Hermitian structures obtained by
Ganchev and Borisov in \cite{Ganchev} are the anti-K\"ahler
structures, studied by the same authors in \cite{Ganchev1}.

In this section we shall give the complete characterization of the
general natural anti-K\"ahlerian structures on the tangent bundle
$(TM,G)$ of a Riemannian manifold $(M,g)$.

To this aim, we shall introduce a $(0,3)$-tensor field $F$ on
$TM$, depending on  the Levi Civita connection $\nabla$ of $G$:
\begin{equation}\label{F}
~~~~~~~~~~~~~F(X,Y,Z)=G((\nabla _X J)Y,Z),\ \forall X,Y,Z\in
\mathcal{T}^1_0(TM),
\end{equation}
with the properties $F(X,Y,Z)=F(X,Z,Y)=F(X,JY,JZ).$

$(TM,G,J)$ is anti-K\"ahlerian if and only if $F=0$. Since $\nabla
$ is torsion free, the condition $F=0$ implies that $J$ is
integrable, hence $(M,g)$ is a space form.

The components of $F$ with respect to the local adapted frame
$\{\delta_i,\partial_j\}_{i,j=\overline{1,n}}$, are denoted by $F$
followed by sequences of $X$ and $Y$, to indicate horizontal or
vertical argument on a certain position. For example, we have
\begin{eqnarray*}
\begin{array}{c}
~~~~~~~~FYYY_{ijk}=F(\partial_i,\partial_j,\partial_k) =G(\nabla
_{\partial_i}J\partial_j,\partial_k)- G(\nabla
_{\partial_i}\partial_j,J\partial_k),
\\
~~~~~~~~FXYX_{ijk}=F(\delta_i,\partial_j,\delta_k) =G(\nabla
_{\delta_i}J\partial_j,\delta_k)- G(\nabla
_{\delta_i}\partial^j,J\delta_k).
\end{array}
\end{eqnarray*}

Using the expression of the Levi-Civita connection of the metric
$G$ (see \cite{OprDruta}), and the relations (\ref{sist4}),
(\ref{defG}), which give the components of $J$ and $G$, we obtain
\begin{eqnarray*}
~~~\begin{array}{c} FYYY_{ijk}= (a_3^\prime c_2-a_2^\prime
c_3+a_2d_3)g_{jk}g_{0i}+(2b_3c_2-2b_2c_3+a_2c_3^\prime
\\
+a_2d_3+2b_3c_2^\prime t)(g_{ik}g_{0j}+ g_{ij}g_{0k})+(b_3^\prime
c_2-b_3c_2^\prime-b_2^\prime c_3+b_2c_3^\prime+a_3^\prime d_2
\\
+ 2b_3d_2-a_2^\prime d_3+a_2d_3^\prime +2b_3^\prime
d_2t-2b_2^\prime d_3t+ 2b_2 d_3^\prime t)g_{0i}g_{0j}g_{0k}.
\end{array}
\end{eqnarray*}

The other components, $FYYX_{ijk}=FYXY_{ikj}$, $FYXX_{ijk}$,
$FXXY_{ikj}=FXYX_{ijk}$, $FXYY_{ijk}$, $FXXX_{ijk},$ depend on the
components $R^h_{kij}$ of the curvature tensor field of the
connection $\dot \nabla$ on $M$. If $(TM,G,J)$ is an anti-K\"ahler
manifold, then $(M,g)$ must have constant sectional curvature c,
i.e. the components of the curvature tensor field of $\dot \nabla$
must be given by
\begin{equation}\label{curbct}
~~~~~~~~~~~~~~~~~~~~~~~~~~~~~~~~R^h_{kij}=c(\delta^h_ig_{kj}-\delta^h_jg_{ki}).
\end{equation}

Taking (\ref{curbct}) into account, and using the relations
(\ref{sist4}) and (\ref{defG}), we obtain
\begin{eqnarray*}
\begin{array}{c}
FYXX_{ijk}=\frac{1}{2}(2b_1c_3 -2b_3c_1 + a_1c_3' + a_1d_3+2b_3cc_2t + 2b_1c_3't-2b_3d_1t \\
+ 2b_1d_3t)(g_{ik}g_{0j}+g_{ij}g_{0k}) + ( a_1'c_3-a_3'c_1+
a_1d_3)g_{jk}g_{0i}+ (b_1'c_3 - b_3'c_1 - b_3cc_2\\
- a_3'd_1 - b_3d_1 + a_1'd_3 + 3b_1d_3 + a_1d_3' - 2b_3'd_1t +
2b_1'd_3t + 2b_1d_3't)g_{0j}g_{0k}g_{0i}.
\end{array}
\end{eqnarray*}

Multiplying the above relation by $y^i,\ y^j,$ and $y^k$,
successively, and taking into account that $g_{0i}y^i=2t,$ it
follows, by using lemma \ref{lema1}, that $FYXX_{ijk}=0$ is
equivalent to the vanishing condition for all its coefficients.
From the coefficient of $g_{jk}g_{0i}$, we obtain that $d_3$ has
the expression:
\begin{equation}\label{vald3}
~~~~~~~~~~~~~~~~~~~~~~~~~~~~~~~~~~~~~d_3=\frac{a^\prime_3c_1-a^\prime_1
c_3}{a_1}.
\end{equation}

The final expression of $FXYY_{ijk}$ is
\begin{eqnarray*}
~~~~~~~~\begin{array}{c} FXYY_{ijk}=-a_2(cc_2 - d_1)g_{jk}g_{0i} +
\frac{1}{2}(a_2c_1' +
a_2cc_2 + a_2d_1 \\
+ 2b_2cc_2t + 2b_3c_3't + 2b_2d_1t - 2b_3d_3t)(g_{ik}g_{0j} +
g_{ji}g_{0k}) \\
 + (b_2c_1' - b_2cc_2 - b_3c_3' + b_2d_1 +
a_2d_1' + b_3d_3 + 2b_2d_1't)g_{0i}g_{0j}g_{0k}.
\end{array}
\end{eqnarray*}

Multiplying the relation $FXYY_{ijk}=0$, by $y^i,\ y^j,\ y^k,$
successively, we obtain by using lemma \ref{lema1}, a simple
expression for $d_1$:
\begin{equation}\label{Vald1}
~~~~~~~~~~~~~~~~~~~~~~~~~~~~~~~~~~~~~~~~~~~d_1=cc_2.
\end{equation}

The final expressions of the components of $F$, become quite long
after replacing the values of the coefficients obtained from the
conditions for $(TM,G,J)$ to be anti-Hermitian.

From the vanishing conditions of the coefficients of
$g_{ik}g_{0j}$ and $g_{ij}g_{0k}$ in the final expression of
$FXXY_{ijk}$, we obtain the values of $c_1'$ and $c_3'$, which
make vanish all the components of $F$. Thus we proved the
following theorem:

\begin{theorem}\label{antikahl}
The tangent bundle $TM$ endowed with the general natural
anti-Hermi-\ tian structure $(G,J)$, is an anti-K\"ahler manifold,
if and only if the almost complex structure $J$ is integrable, the
coefficients $d_1,\ d_3$ from the definition of $G$ are given by
$(\ref{Vald1})$ and $(\ref{vald3})$, respectively, and $c_1',\
c_3',$ have the values \large\begin{eqnarray*}
~~~~\begin{array}{c} c_1'=-2ca_1^2\frac{2 c_3 [ a_1 (a_3 + a_3' t)
- 2 a_1' a_3 t] - c_1 (1 + a_3^2)}{a_1^4 + 4 a_1^2 (a_3^2-1) c t +
4 (1 + a_3^2)^2 c^2 t^2}
\\
-2c(1 + a_3^2) \frac{c c_1 t [2 a_1 (1 + a_3^2 + 4 a_3 a_3' t)-4
a_1'(1 + a_3^2)  t ] +
    2 a_1^2 t (a_1' c_1 - 2 a_3' c c_3 t)}{a_1^5 + 4 a_1^3 (a_3^2-1) c t + 4 a_1 (1 + a_3^2)^2
c^2 t^2},
\\\\
c_3'=2 c a_1  \frac{c_3 \{a_1 +
      a_3 [4 t (a_1' a_3  - a_1 a_3')-3 a_1 a_3]\} +
         2 a_3 c_1 (1 + a_3^2 + 2 a_3 a_3' t)}{a_1^4 + 4 a_1^2 (a_3^2-1) c t +
  4 (1 + a_3^2)^2 c^2 t^2}
\\
   +\frac{(a_3' c_1 - a_1' c_3) [a_1^4 -
      4 (1 + a_3^2)^2 c^2 t^2] -
   4 a_1 (1 + a_3^2) c [2 a_1' a_3 c_1 + (1 + a_3^2) c c_3] t}{a_1^5 + 4 a_1^3 (a_3^2-1) c t +
  4 a_1 (1 +
     a_3^2)^2 c^2 t^2}.
\end{array}
\end{eqnarray*}

\normalsize If the sectional curvature of the base manifold $M$ is
a positive constant, c, the anti-K\"ahlerian structure is defined
on the whole tangent bundle $TM$, and if c is strictly negative,
the anti-K\"ahlerian structure is defined only in the condition
$$
(a_1^2 - 2ct)^2  \neq -4a_3^2ct(a_1^2+ 2ct + a_3^2ct).
$$
\end{theorem}

In a similar way, we may prove the following theorem, which
characterize the anti-K\"ahlerian tangent bundles of natural
diagonal lift type.

\begin{theorem}
The almost anti-Hermitian manifold $(TM,G,J)$ of diagonal lift
type is an anti-K\"ahler manifold if and only if

1) The natural diagonal almost complex structure $J$ is integrable
$($i. e. the base manifold $M$ is of constant sectional curvature,
c, and $b_1=\frac{a_1a_1'- c}{a_1- 2ta_1'}$$)$, and
\begin{eqnarray}\label{d1c1p}
~~~~~~~~~~~~~~\begin{array}{c} c_1'=\frac{2cc_1(a_1 -
2a_1't)}{a_1(a_1^2 - 2ct)},\ d_1 = c c_2 = -c \frac{c_1}{a_1^2}, \
a_1^2\neq 2ct.
\end{array}
\end{eqnarray}

2) The base manifold has constant sectional curvature $c$, and the
essential coefficients have the expressions:
$$
a_1=\sqrt{B+2ct},\ b_1=0,\ c_1=A(B+2ct),\ d_1=-cA,
$$
where  $A$ is a nonzero real constant and $B$ is a positive
constant.

If $c$ is positive, the anti-K\"ahler structure from the second
case is defined on the all $TM$, and if $c$ is strictly negative,
the structure is defined only on the tube $ t<-\frac{B}{2c} $
around the null section in $TM.$
\end{theorem}

\section{Conformally anti-K\"ahler structures of general na-\ tural lift type on the tangent bundle}

The conformally anti-K\"ahler structures or $\omega_1-$structures
on the tangent bundle of a Riemannian manifold $(M,g)$ are
characterized in \cite{Ganchev} by the relation
\begin{eqnarray}\label{conf}
~~~~~~~~~~~~~~\begin{array}{c}
2nF(X,Y,Z)=G(X,JY)\Phi(JZ)+G(X,JZ)\Phi(JY)\\
+G(X,Y)\Phi(Z)+G(X,Z)\Phi(Y),\ \forall X,Y,Z\in
\mathcal{T}^1_0(TM),
\end{array}
\end{eqnarray}
where $F$ is the usual tensor field of type (0,3), given by
(\ref{F}), and $\Phi$ is a 1-form associated with $F$, which in
the case of the tangent bundles with general natural almost
anti-Hermitian structures, is given by
\begin{eqnarray}\label{Phi}
~~~~~~~~~~~~~~~~~~~~~\begin{array}{c}
\Phi(X)=H_{(1)}^{ij}F(\delta_i,\delta_j,X)+H_{(3)}^{ij}F(\delta_i,\partial_j,X)\\
+H_{(3)}^{ij}F(\partial_i,\delta_j,X)
+H_{(2)}^{ij}F(\partial_i,\partial_j,X),
\end{array}
\end{eqnarray}
$\forall X\in \mathcal{T}^1_0(TM), \forall i,j=\overline{1,n},$
$H$ being the inverse matrix of $G$ (see \cite{OprDruta}).

\vskip2mm In the following theorem we shall characterize the
tangent bundles endowed with conformally anti-K\"ahler structures
of general natural lift type.

\begin{theorem}\label{confantikahl}
The almost anti-Hermitian manifold $(TM,G,J)$, is general natural
conformally anti-K\"ahler, if and only if the almost complex
structure $J$ is integrable $($the base manifold is of constant
sectional curvature c, and $b_1,$ $b_2,$ $b_3$ have the
expressions $(\ref{integrab})$$)$, and $d_1,\ d_3$ are of the
forms
\begin{eqnarray}\label{d1d3}
~~~~~~~~~~~~~~~~~~~~\begin{array}{c}
d_1=c\frac{2a_1a_3c_3-c_1(1+a_3^2)}{a_1^2},\ d_3=\frac{a_3'c_1 -
a_1'c_3}{a_1}.
\end{array}
\end{eqnarray}
\end{theorem}

\emph{Proof:} The vertical and horizontal components of the 1-form
$\Phi$ are
\begin{eqnarray*}
\Phi\partial_k=H_{(1)}^{ij}FXXY_{ijk}+H_{(3)}^{ij}FXYY_{ijk}+H_{(3)}^{ij}FYXY_{ijk}+H^{ij}_{(2)}FYYY_{ijk},
\\
\Phi\delta_k=H_{(1)}^{ij}FXXX_{ijk}+H_{(3)}^{ij}FXYX_{ijk}+H_{(3)}^{ij}FYXX_{ijk}+H^{ij}_{(2)}FYYX_{ijk}.
\end{eqnarray*}

If we replace $X,Y,Z$ by $\partial_i,\partial_j,\delta_k,$ we
obtain from (\ref{conf}) the relation
\begin{eqnarray}\label{confFYYX} ~~~\begin{array}{c}
2nFYYX_{ijk}=[G^{(2)}_{il}(J_1)_k^l-G^{(3)}_{il}(J_3)^l_k][(J_3)^m_j\Phi\partial_m-(J_2)^m_j\Phi\delta_m]
\\
+[G^{(2)}_{il}(J_3)_j^l-G^{(3)}_{il}(J_2)^l_j][(J_1)_k^m\Phi\partial_m-(J_3)^m_k\Phi\delta_m]+G^{(2)}_{ij}\Phi\delta_k+G^{(3)}_{ik}\Phi\partial_j,
\end{array}
\end{eqnarray}
\normalsize which becomes quite complicated after replacing the
components of $J$, $G$ , and $\Phi$. Multiplying the expression by
$g^{il}g^{jm}$, and taking into account that the curvature of the
base manifold do not depend on the tangential coordinates, we
compute the derivative of the final expression with respect to
$y^k$. Taking the value in $y=0$, we obtain an expression which
depends on the components of the curvature of the base manifold,
and on those of the Ricci tensor, which leads to the Einstein
condition for the base manifold. Using this condition and the
first Bianchi identity, we obtain that the base manifold is a
space form, i.e. its curvature has the expression (\ref{curbct}),
where c depends on the values in $t=0$ of $a_1,\ a_1',\ a_3,\
a_3',\ c_1,\ c_3,\ b_1,\ b_3.$

Using the theorems \ref{antiherm}, \ref{th4}, and the relation
(\ref{curbct}), the condition (\ref{confFYYX}) becomes
~~~~~~~~\begin{eqnarray}\label{nouFYYX} Ag_{ij}g_{0k}+
\frac{2na_3}{a_1}(a_1'c_3 - a_3'c_1 + a_1d_3)g_{jk}g_{0i} +
Bg_{ik}g_{0j}+Cg_{0k}g_{0i}g_{0j}=0,
\end{eqnarray}
where $A,\ B,\ C$ are some functions depending on the coefficients
of $J$ and $G$.

Multiplying (\ref{nouFYYX}) by $y^k,\ y^i,\ y^j$ successively, and
using lemma \ref{lema1}, we obtain that (\ref{nouFYYX}) is
equivalent to the vanishing condition for all the coefficients
involved. Thus, $d_3$ must have the expression from (\ref{d1d3}),
given in the theorem.

Writing the relation (\ref{conf}) for $\delta_i,\
\partial_j,\ \partial_k$ instead of $X,Y,Z,$ respectively, we
get the value of $d_1$, in a similar way as $d_3$ was obtained.

Replacing (\ref{d1d3}) into the vanishing condition for the
coefficient $A$, we get
$$
b_3=\frac{a_1a_3' [a_1^2- 2ct (1+ 3 a_3^2)]  + 2 a_3 c (a_1- 2a_1'
t) (1+ a_3^2)}{(a_1 - 2 a_1't) [a_1^2 - 2ct (1 + a_3^2)]
 - 8a_1a_3a_3'ct^2}.
$$

Similarly, from (\ref{conf}), written for $FYXX_{ijk}$, we have
$$
b_1=\frac{a_1a_1' [a_1^2- 2ct(1+ a_3^2)]  - a_1^2c [1 - a_3 (3a_3+
4a_3't)] + 2c^2t(1 + a_3^2)^2}{(a_1 - 2 a_1't) [a_1^2 - 2ct (1 +
a_3^2)]
 - 8a_1a_3a_3'ct^2}.
$$

\normalsize We may easily verify, that the above values and of
$b_1$ and $b_3$, and the value of $b_2$ from (\ref{inlocuire})
coincide with those given in the integrability conditions
(\ref{integrab}), in which we replace $a_2$ from
(\ref{inlocuire}). Thus the theorem is proved.
\qquad\qquad\qquad\qquad\qquad\qquad\qquad$\square$

\bf{Remark} 4.1 \rm Particularizing the result from the theorem
\ref{confantikahl} to the diagonal case, we obtain the
characterization of the conformally anti-K\"ahler structures of
natural diagonal lift type on $TM$, given in [24, theorem 14].

\section{General natural complex anti-Hermitian structures on $TM$}
The complex anti-Hermitian structures, or
$\omega_1\oplus\omega_2-$ structures on the tangent bundle of a
Riemannian manifold $(M,g)$, are characterized by the condition
\begin{equation}\label{compl}
~~~~~~\begin{array}{c} F(X,Y,JZ)+F(Y,Z,JX)+F(Z,X,JY)=0, \forall
X,Y,Z\in \mathcal{T}^1_0(TM),
\end{array}
\end{equation}
where $F$ is the usual tensor field of type (0,3), given by
(\ref{F}), and $J$ is the almost complex structure of general
natural lifted type, defined by $(\ref{sist4})$.

In this section we shall prove the characterization theorem for
the complex anti-Hermitian manifolds $(TM,G,J)$ of general natural
lift type.

\begin{theorem}\label{thcaractcomplex}
The general natural anti-Hermitian manifold $(TM,G,J)$, is complex
anti-Hermitian, if and only if the almost complex structure $J$ is
integrable.
\end{theorem}

\emph{Proof:} Since the curvature of the base manifold satisfy the
Bianchi identity, the relation (\ref{compl}) is verified when we
replace $X,\ Y,\ Z$ by $\partial_i,\
\partial_j,\ \partial_k$, or by $\delta_i,\ \delta_j,\ \delta_k,$
respectively. Thus, the tangent bundle endowed with general
natural complex anti-Hermitian structure will be characterized by
the following essential relations only:
\begin{equation}\label{compl1}
~~~~~~~~~~~~F(\partial_i,\partial_j,J\delta_k)+F(\partial_j,\delta_k,J\partial_i)+F(\delta_k,\partial_i,J\partial_j)=0,
\end{equation}
\begin{equation}\label{compl2}
~~~~~~~~~~~~F(\partial_i,\delta_j,J\delta_k)+F(\delta_j,\delta_k,J\partial_i)+F(\delta_k,\partial_i,J\delta_j)=0.
\end{equation}

The value in $y=0$ of the derivative of (\ref{compl1}) with
respect to the tangential coordinates $y^k$, leads to the
condition for the base manifold to be a space form. Then, taking
(\ref{curbct}) into account, we may write the relations
(\ref{compl1}) and (\ref{compl2}) in the forms $f_1
(g_{jk}g_{0i}-g_{ik}g_{0j})=0,\ f_2(g_{ik}g_{0j}-g_{ij}g_{0k})=0,$
where $f_1$ and $f_2$ are rational functions depending on $a_1,\
b_1,\ c_1,\ a_3,\ b_3,$ $c_3,$ $a_1',\ a_3'$, on the constant
sectional curvature $c$, and on the energy density, $t$. By using
lemma \ref{lema1}, we have that $f_1$ and $f_2$ must vanish.
Taking into account these vanishing conditions, and considering
the expression (\ref{inlocuire}) of $b_2$, we obtain the
integrability conditions from theorem $\ref{th3}$.
\qquad\qquad\qquad\qquad\quad $\square$

\bf{Remark} 5.1 \rm If in the theorem \ref{thcaractcomplex} we
consider the anti-Hermitian structure $(G,J)$ of diagonal lift
type on the tangent bundle $TM$, we obtain the result from [24,
theorem 9], which characterize the natural diagonal complex
anti-Hermitian tangent bundles.

\section{Tangent bundles with general natural quasi-anti-K\"ahlerian structures}

In this section we shall obtain the conditions under which the
tangent bundle $TM,$ with a general natural anti-Hermitian
structure $(G,J)$ is quasi-anti-K\"ahler.

The quasi-anti-K\"ahlerian manifolds, or $\omega_3-$ manifolds,
are characterized in \cite{Ganchev} by the following vanishing
condition:
\begin{equation}\label{quasianti}
~~~~~~~~~~~F(X,Y,Z)+F(Y,Z,X)+F(Z,X,Y)=0, \forall X,Y,Z\in
\mathcal{T}^1_0(TM).
\end{equation}
where $F$ is the (0,3)-tensor field, given by (\ref{F}).

If in (\ref{quasianti}) we take $\partial_i,\ \delta_j,\
\delta_k$, instead of $X,Y,Z,$ respectively, we have
$$
FYXX_{ijk} + FXXY_{jki} + FXYX_{kij}=0,
$$
relation which depends on the curvature of the base manifold.
Differentiating the final form of the above relation, with respect
to the tangential coordinates $y^h$, and taking the value in
$y=0,$ we obtain by standard calculation that the base manifold is
a space form, i.e $R^h_{kij}$ has the form (\ref{curbct}). Then,
the condition
$$
FXXX_{ijk}+FXXX_{jki}+FXXX_{kij}=0,
$$
becomes of the form
\begin{eqnarray}\label{FQQQquasi}
~~~~~~~~~\begin{array}{c} (a_1 c_1' + 2 a_1 d_1 + 2 b_1 c_1' t + 4
b_3 c c_3 t) (g_{jk}g_{0i} + g_{ik} g_{0j} +  g_{ij}g_{0k})
\\
-3 (2 b_3 c c_3 - 2 b_1 d_1 - a_1 d_1' - 2 b_1 d_1' t)
g_{0i}g_{0j}g_{0k}=0.
\end{array}
\end{eqnarray}
Using lemma \ref{lema1}, we obtain from (\ref{FQQQquasi}) the
expressions of $c_1'$ and $d_1'$:
\begin{eqnarray}\label{c1pquasi}
~~~~~~~~~~~~~~~~~~~~\begin{array}{c} c_1'=-\frac{2(a_1d_1 +
2b_3cc_3t)}{a_1 + 2b_1t},\ d_1'=\frac{2(b_3cc_3 - b_1d_1)}{a_1 +
2b_1t}.
\end{array}
\end{eqnarray}

Replacing (\ref{c1pquasi}) into the other components of the sum
(\ref{quasianti}), we obtain by similar computations, the values
of $a_1',\ a_3',\ c_3',\ d_3'$, and we may formulate:

\begin{theorem}\label{quasiantiK}
The tangent bundle $TM$ endowed with the almost anti-Hermitian
structure $(G,J)$ of general natural lift type, is a
quasi-anti-K\"ahler manifold, if and only if the base manifold is
of constant sectional curvature $c$, the first order derivatives
of $c_1$ and $d_1$ have the expressions $(\ref{c1pquasi})$, and
the other coefficients verify \large
\begin{eqnarray*}
\begin{array}{c}
a_1'=
  a_1\frac{b_1 c_1 (1 + a_3^2 + 2 a_3 a_3' t) +
     2 (1 + a_3^2) c c_3 (a_3 - b_3 t)}{[c_1 (1 + a_3^2) - a_1 a_3 c_3] (a_1 +
    2 b_1 t)}-\frac{(1 + a_3^2)^2 c c_1 + a_1^3 [(a_3' - b_3) c_3 + a_3
c_3']}{[c_1 (1 + a_3^2) - a_1 a_3 c_3] (a_1 +
    2 b_1 t)}
    \\
  -a_1^2\frac{2 d_1 + a_3 (2 b_1 c_3 + 2 a_3 d_1-a_3' c_1) +
     2 b_1 (a_3' c_3 + a_3 c_3') t}{[c_1 (1 + a_3^2) - a_1 a_3 c_3] (a_1 +
    2 b_1 t)},
\\\\
a_3'= -
   a_1 (1 + a_3^2) \frac{b_3 (c_1^2 + 2 c c_3^2 t + 4 c_1 d_1 t) +
      2 c_1 [a_3 d_1 - b_1 (c_3 + 2 d_3 t)]}{[(1 + a_3^2) c_1^2 -
     a_1 c_3 (2 a_3  c_1  + a_1 c_3)] (a_1 + 2 b_1 t)}
    \\
    + 2 a_1^2 \frac{c_1 d_3 - c_3 d_1 + a_3^2 (c_3 d_1 + c_1 d_3) -
      a_3 c_3 [b_1 (c_3 + 2 d_3 t)-2 b_3 d_1 t]}{[(1 + a_3^2) c_1^2 -
     a_1 c_3 (2 a_3  c_1  + a_1 c_3)] (a_1 + 2 b_1 t)}
     \\
   +\frac{a_1^4 c_3 (b_3 c_3 - 2 a_3 d_3) - 2 (1 + a_3^2) b_3 c
c_1 t [c_1 (1 + a_3^2) +2 a_1 a_3 c_3]}{a_1 [(1 + a_3^2) c_1^2 -
     a_1 c_3 (2 a_3  c_1  + a_1 c_3)] (a_1 + 2 b_1 t)} ,
\\\\
c_3'=\frac{2 (1 + a_3^2) b_3 c c_1 t}{a_1^2 (a_1 + 2 b_1
t)}+\frac{a_1 (b_3 c_1 - b_1 c_3 - 2 a_3 d_1) +
 c c_3 (1 + a_3^2 - 4 a_3 b_3 t)}{a_1(a_1 + 2 b_1 t)},
\end{array}
\end{eqnarray*}
\normalsize and $d_3'$ has a longer and more complicated
expression of the same type.

The result is true only when the non-vanishing conditions for the
denominators of $a_1'$ and $a_3'$ is satisfied.
\end{theorem}

\bf{Remark} 6.1 \rm Considering the case of the natural diagonal
quasi-anti-K\"ahler structures on the tangent bundle, from the
above theorem we obtain the characterization given in [24, theorem
8].

\section{Tangent bundles with semi-anti-K\"ahler structures of general natural lift type}

The semi-anti-K\"ahler manifolds, or $\omega_2\oplus \omega_3-$
manifolds, are defined in \cite{Ganchev} as being the manifolds
for which the 1-form $\Phi$, given by (\ref{Phi}), vanishes:
\begin{equation}
~~~~~~~~~~~~~~~~~~~~~~~~~~~~~~~~~~~~~~~~~~~\Phi=0.
\end{equation}

We shall prove the following proposition, in which we give some
necessary conditions for the anti-Hermitian manifold $(TM,G,J)$ of
general natural lift type to be semi anti-K\"ahler.

\begin{proposition}\label{semic3p}
If the general natural almost anti-Hermitian structure $(G,J)$ on
the tangent bundle $TM$ is semi-anti-K\"ahler, then the base
manifold is Einstein, and the first order derivative of the
coefficient $c_3$ from the definition of $G$ is given by \large
\begin{eqnarray*}
\begin{array}{c}
c_3'=\frac{2 a_1' (1 + a_3^2) b_3 c_1^2 t-a_1^2c_3 \{c_1'(1 -
a_3^2) - 2 a_3' b_3 c_1 t + a_3 [a_3' c_1 + a_1' c_3 - 2 b_3 (c_1
+ c_1' t)]-a_1c_3 (a_3' - b_3)\}}{a_1^2 [a_1 c_3 (a_3 +
        2 b_3 t) - c_1 (1 + a_3^2 + 2 a_3 b_3 t)]}
\\
-\frac{c_1}{a_1}\frac{(1 + a_3^2) (b_3 c_1 + a_3 c_1' - a_1' c_3)+
2 b_3 (c_1' + a_3 (a_3' c_1 + a_3 c_1' + a_1' c_3)) t}{a_1 c_3
(a_3 + 2 b_3 t) - c_1 (1 + a_3^2 + 2 a_3 b_3 t)},
\end{array}
\end{eqnarray*}
\normalsize when $a_1 c_3 (a_3 + 2 b_3 t) \neq c_1 (1 + a_3^2 + 2
a_3 b_3 t)$.
\end{proposition}

\emph{Proof:} The tangent bundle $TM$, endowed with the general
natural almost anti-Hermitian structure $(G,J)$ is a
semi-anti-K\"ahler manifold if and only if the tensor $\Phi$
vanishes, i. e. if and only if $\Phi\delta_k=\Phi\partial_k=0$.
The final expression of $\Phi\delta_k,$ is
\begin{eqnarray}\label{Phidelta}
~~~~~~~~~~~~~~~~~~~~~~~~~~\begin{array}{c} \Phi\delta_k = f(t)
g_{0k}-\frac{1 + a_3^2}{a_1}y^h Ric_{hk},
\end{array}
\end{eqnarray}
where $f$ is a rational function depending on the coefficients of
$G$ and $J$, on the energy density, $t$, and on the dimension $n$
of $M$. Since the Ricci tensor of the base manifold, \emph{Ric} do
not depend on the tangential coordinates, the value in $y=0$ for
the derivative of $\Phi\delta_k$ with respect to $y^h$ leads to
the Einstein condition for the base manifold:
$$
Ric_{hk}=\rho g_{hk},
$$
where $\rho=f(0)\frac{a_1(0)}{1+a_3(0)^2}.$

The vertical component, $\Phi\partial_k$ has the form
$(u(t)+nv(t))g_{0k}=0,$ $u$ and $v$ being two rational functions
depending on the same parameters as $f$.

Since $\Phi\partial_k$ must vanish for every dimension $n$ of the
base manifold, we have that $u=v=0,$ and solving the equation
$v=0$ with respect to $c_3'$, we obtain the expression from the
proposition.~\qquad\qquad\qquad\qquad\qquad\qquad\qquad\qquad
\qquad\qquad\qquad\qquad\quad\quad\quad$\square$

\bf{Remark} 7.1 \rm The sufficient conditions for the general
natural anti-Hermitian ma-\ nifold $(TM,G,J)$ to be
semi-anti-K\"ahler, are given by the Einstein property of the base
manifold, by the value of $c_3'$, and by other two more
complicated relations between the coefficients of $G$ and $J$,
which may not be presented here.

\bf{Remark} 7.2 \rm If the base manifold is Ricci flat and the
general natural almost anti-Hermitian manifold $(TM,G,J)$ is
semi-anti-K\"ahler, then $c_3'$ has the expression from
proposition \ref{semic3p}, and when $c_1 (1 + a_3^2) \neq a_1 a_3
c_3$, $a_1'$ has the form: \large\begin{eqnarray*}
\begin{array}{c}
a_1'=a_1^2 \frac{a_3 [ c_1 (a_3' - b_3) - b_1 c_3] + c_1' - 2 a_3'
b_1 c_3 t+ c_3(b_3-a_3')}{[c_1 (1 + a_3^2) -
   a_1 a_3 c_3] (a_1 + 2 b_1 t)}+ b_1 \frac{2 c_1' t + c_1 (1 + a_3^2 + 2
a_3 a_3' t)}{[c_1 (1 + a_3^2) -
   a_1 a_3 c_3] (a_1 + 2 b_1 t)}.
\end{array}
\end{eqnarray*}

\normalsize \bf{Remark} 7.3 \rm For the natural diagonal
semi-anti-K\"ahler structure $(G,J)$ on $TM$, the expressions
which we have to study become simpler, and the relation between
the coefficients of $G$ and $J$ is that obtained in [24, theorem
10], which has an interesting consequence, furnishing a simple
example of semi-anti-K\"ahler structure, presented in the same
paper.

\section{Special complex anti-Hermitian structures of ge-\ neral natural
lift type on the tangent bundle}

The complex anti-Hermitian manifolds which are at the same time
semi-anti-K\"ahlerian manifolds are called in \cite{Ganchev},
special complex anti-Hermitian manifolds, or $\omega_2 -$
manifolds.

The special complex anti-Hermitian manifolds $(TM,G,J)$, of
general natural lift type are characterized by the conditions
$$
\Phi=0,\ F(X,Y,JZ)+F(Y,Z,JX)+F(Z,X,JY)=0, \forall X,Y,Z\in
\mathcal{T}^1_0(TM).
$$
where $F$ is the usual tensor field of type (0,3), defined by
(\ref{F}), $\Phi$ is the 1-form associated with $F$, given by
(\ref{Phi}), and $J$ is the almost complex structure of general
natural lifted type with the expression $(\ref{sist4})$.

In this section we shall find some necessary conditions under
which the tangent bundle $TM$ endowed with a general natural
anti-Hermitian structure is special complex anti-Hermitian.

\begin{proposition}\label{specialc1p}
If the anti-Hermitian manifold $(TM,G,J)$ of general natural lift
type is special complex anti-Hermitian, then the almost complex
structure $J$ is integrable $($i.e. the base manifold is a space
form, and $b_1,\ b_2,\ b_3$ have the values $(\ref{integrab})$$)$,
$c_3'$ has the expression from proposition $\ref{semic3p}$, and
$c_1'$ is given by \large\begin{eqnarray*} ~~~~~~
\begin{array}{c}
c_1'=2 c \frac{a_1^3 [2 a_1 c_3 (a_3 + a_3' t) -
   c_1 (1 + a_3^2) - 4 a_1' a_3 c_3 t]-4 a_1' (1 + a_3^2)^2 c c_1 t^2 }{
 a_1^5 + 4 a_1^3 (a_3^2-1) c t + 4 a_1 (1 + a_3^2)^2 c^2 t^2}
\\
-4 c t\frac{(1 + a_3^2) [c c_1 (1 + a_3^2 + 4 a_3 a_3' t)+ a_1
(a_1' c_1 - 2 a_3' c c_3 t)]}{a_1^4 + 4 a_1^2 (a_3^2-1) c t + 4 (1
+ a_3^2)^2 c^2 t^2}.
\end{array}
\end{eqnarray*}

\normalsize For $c\geq 0$ the result is always valid, and if $c>0$
the coefficients must verify the relation
$[a_1+2ct(1+a_3^2)]^2\neq8a_1ct.$
\end{proposition}

\emph{Proof:} Since the special complex anti-Hermitian manifold
$(TM,G,J)$ is complex anti-Hermitian and semi-anti-K\"ahler, we
have from theorem \ref{thcaractcomplex} that the base manifold $M$
has constant sectional curvature $c$, and from proposition
\ref{semic3p} we obtain that $M$ is Einstein, the constant $\rho$
being in this case equal to $c(n-1)$. Thus we obtain that
$\Phi\delta_k=(\alpha(t)+\beta(t) n)g_{0k}=0,$ where $\alpha$ and
$\beta$ are rational functions depending on the energy density and
on the coefficients of the metric $G$ and of the almost complex
structure $J$. After replacing the value of $c_3'$ and the
expressions (\ref{integrab}) and (\ref{inlocuire}) into the
relation $\beta=0$, we obtain the value of $c_1'$ from the
proposition.\qquad\qquad\qquad\qquad$\square$

\bf{Remark} 8.1 \rm The sufficient conditions for the general
natural anti-Hermitian ma-\ nifold $(TM,G$, $J)$ to be special
complex anti-Hermitian, are given by the integrability of the
almost complex structure $J$, by the expressions  of $c_1'$ and
$c_3'$ from the propositions $\ref{specialc1p}$, $\ref{semic3p}$,
and by two complicated relations between the coefficients $G$ and
$J$, which may not be presented here.

\bf {Remark 8.2} \rm When the almost anti-Hermitian manifold
$(TM,G,J)$ is of natural diagonal lift type, the sufficient
conditions for $(TM,G,J)$ to be special complex anti-Hermitian
become simpler, the complete characterization of this structures
being given in [24, theorem 12].

\section{$\omega_1\oplus\omega_3$ structures on $TM$}

An almost anti-Hermitian manifold $(TM,G,J)$ is an
$\omega_1\oplus\omega_3-$ manifold, if it satisfies the relation
\small\begin{eqnarray}\label{omega} ~~~~~~~~~~~~~~
\begin{array}{c}
n[F(X,Y,Z)+F(Y,Z,X)+F(Z,X,Y)]=\\
=G(X,Y)\Phi(Z)+G(Z,X)\Phi(Y)+G(Y,Z)\Phi(X)\\
+G(X,JY)\Phi(JZ)+G(Y,JZ)\Phi(JX)+G(Z,JX)\Phi(JY),
\end{array}
\end{eqnarray}
\normalsize where $X,Y,Z$ are any vector fields on $TM$, $F$ is
the usual tensor field of type (0,3), with the expression
(\ref{F}), $\Phi$ is the associated 1-form, defined by
(\ref{Phi}), $G$ and $J$ are the general natural semi-Riemannian
metric and almost complex structutre, given by $(\ref{defG})$ and
$(\ref{sist4})$, respectively.

In the following proposition we shall present some necessary
conditions for the almost anti-Hermitian manifold $(TM,G,J)$ of
general natural lift type to be an $\omega_1\oplus\omega_3-$
manifold.

\begin{proposition}\label{omega1+omega3}
If the almost anti-Hermitian manifold $(TM,G,J)$ is an
$\omega_1\oplus\omega_3-$ manifold of general natural lift type,
then the base manifold has constant sectional curvature $c$,
\large
\begin{eqnarray*}
~~~\begin{array}{c} a_1'=\frac{a_1^3 (b_3-a_3') c_3 - (1 + a_3^2)
c c_1 (1 + a_3^2 + 2 a_3 b_3 t)} {[c_1 (1 +
   a_3^2) - a_1 a_3 c_3] (a_1 + 2 b_1 t)} +
 a_1^2 \frac{a_3 [c_1 (a_3'  - b_3) - b_1 c_3] - 2 (d_1 + a_3' b_1 c_3 t)}{[c_1 (1 +
   a_3^2) - a_1 a_3 c_3] (a_1 + 2 b_1 t)}
\\
   +
a_1 \frac{b_1 c_1 (1 + a_3^2 + 2 a_3 a_3' t) +
    c c_3 [a_3 (1 + a_3^2) - 2 b_3 t (1 - a_3^2)]}{[c_1 (1 +
   a_3^2) - a_1 a_3 c_3] (a_1 + 2 b_1 t)},
\\\\
a_3'=\frac{a_1^4 c_3 (b_3 c_3 - 2 a_3 d_3) - 2 (1 + a_3^2) b_3 c
c_1 t[(1 + a_3^2)c_1 + 4 a_1 a_3  c_3]}{a_1 [(1 + a_3^2) c_1^2 -
      a_1 c_3 (2 a_3  c_1+ a_1 c_3)] (a_1 + 2 b_1 t)}
\\
    -
   a_1^2 (1 + a_3^2) \frac{b_3 (c_1^2 + 2 c c_3^2 t + 4 c_1 d_1 t) +
      2 c_1 [a_3 d_1 - b_1 (c_3 + 2 d_3 t)]}{a_1 [(1 + a_3^2) c_1^2 -
      a_1 c_3 (2 a_3  c_1+ a_1 c_3)] (a_1 + 2 b_1 t)}
\\
      +
   2 a_1^3 \frac{ c_1 d_3 - c_3 d_1 + a_3^2 (c_3 d_1 + c_1 d_3) -
      a_3 c_3 [b_1 (c_3 + 2 d_3 t)-2 b_3 d_1 t]}{a_1 [(1 + a_3^2) c_1^2 -
      a_1 c_3 (2 a_3  c_1 + a_1 c_3)] (a_1 + 2 b_1 t)},
\end{array}
\end{eqnarray*}
\normalsize and other more complicated relations between the
coefficients are satisfied, for example the first order
derivatives of $d_1$ and $d_3$ may be expressed as rational
functions of the other coefficients, their derivatives, c and t.

The proposition is true only when the denominators do not vanish.
\end{proposition}

\emph{Proof:} If we write the relation (\ref{omega}) for the
vector fields $\partial_i,\partial_j,\partial_k,$ we obtain an
expression depending on the tangential coordinates $y^h$, on the
components $g_{ij}$ of the metric $g$, on the entries $g^{ij}$ of
the inverse matrix of $g$, and on the Ricci tensor of the base
manifold. From the value in $y=0$ of the derivative with respect
to $y^h$ of this expression, we obtain, after the multiplication
by $g^{ik}g^{jl}g^{hm}$, that the base manifold is Einstein.

Writting (\ref{omega}) for the vector fields $\partial_i,\
\partial_j,\ \delta_k$, we get, by taking into account the Einstein condition
for the base manifold, that $M$ has constant sectional curvature,
c.

From $(\ref{omega})$ written for $\delta_i,\delta_j,\delta_k$, we
obtain an expression of the form
$$
(u_1+u_2n)(g_{jk}g_{0i}+g_{ik}g_{0j}+g_{ij}g_{0k})+(v_1+v_2n)g_{0i}g_{0j}g_{0k}=0,
$$
where $u_1,\ u_2,\ v_1,\ v_2$ are rational functions depending on
the coefficients of $G$ and $J$. Since the expression must vanish
for every dimension $n$ of $M$, $u_1,\ u_2,\ v_1,\ v_2$ must be
zero. From this vanishing conditions and from the similar
relations obtained from the other components in (\ref{omega}), we
obtain the expressions from the proposition.\qquad $\square$

\bf {Remark} 9.1 \rm The sufficient conditions under which the
anti-Hermitian manifold $(TM,G,J)$ is a of general natural
$\omega_1\oplus \omega_3-$ manifold, are given by the condition
for the base manifold to have constant sectional curvature, c, by
the expressions of $a_1',$ $a_3',$ from proposition
\ref{omega1+omega3}, and by other more complicated relations
between the coefficients of $G$ and $J$, such as some quite long
expressions of $d_1'$ and $d_3'$.

\bf{Remark} 9.2 \rm In the diagonal case it is easier to present
all the necessary and sufficient conditions for $(TM,G,J)$ to be
an $\omega_1\oplus \omega_3-$ manifold (see [24, theorem 15] and
the consequence, which furnishes a simple example of
$\omega_1\oplus \omega_3-$ manifold).

\textbf{Acknowledgements.} The author expresses her gratitude to
Professor Oproiu for the technics learned and for the useful
discussions throughout this work.

Simona-Luiza Dru\c t\u a

Faculty of Mathematics

"AL. I. Cuza" University of Ia\c si

Bd. Carol I, No. 11

RO-700506 Ia\c si

e-mail: simonadruta@yahoo.com

\end{document}